\documentclass[12pt]{amsart}
\usepackage[centertags]{amsmath}
\usepackage{amsfonts}
\usepackage{amssymb}
\usepackage{amscd}
\usepackage{amsthm}
\usepackage{fancyhdr}
\usepackage[dvips]{graphicx}

\def\ra{{\sf ran}}
\def\je{{\sf ker}}
\def\I{\mathbf{I}}
\def\mA{\mathbf{A}}

\def\mR{\mathbf{R}}

\def\mT{\mathbf{T}}

\def\mH{\mathbf{H}}
\def\x{\mathcal{X}}
\def\q{\mathcal{Q}}

\def\d{\mathcal{D}}

\def\H{\mathcal{H}}

\def\R{ \mathbb{R}}
\def\C{ \mathbb{C}}
\def\N{ \mathbb{N}}

\newtheorem{theorem}{Theorem}[section]
\newtheorem{corollary}[theorem]{Corollary}
\newtheorem{lemma}[theorem]{Lemma}

\theoremstyle{definition}
\newtheorem{definition}[theorem]{Definition}
\newtheorem{remark}[theorem]{Remark}
\newtheorem{example}[theorem]{Example}

\numberwithin{equation}{section}

\newcommand{\norm}[1]{\Vert #1 \Vert}

\newcommand{\sk}[1]{\left( #1 \right)}

\newcommand{\abs}[1]{\left| #1 \right|}

\def\sinbf{\mathsf{sin}}

\begin{document}
\title{On eigenvalue and eigenvector estimates for nonnegative definite operators}

\date{\today}

\author{Luka Grubi\v{s}i\'{c}
}
\thanks{This work is a part of author's PhD thesis which was written
under the supervision of Prof. Dr. Kre\v{s}imir Veseli\'{c}, Hagen and is submitted
to the Fachbereich Mathematik der Fernuniversit\"{a}t in Hagen in partial
fulfillment of the requirements for the degree Dr. rer. nat.}
\address{
Luka Grubi\v{s}i\'{c}\\
FernUniversit\" at in Hagen\\
LG Mathematische Physik\\
Feithstr.~140\\
D-58084~Hagen, Germany
}
\email{luka.grubisic@fernuni-hagen.de}

\keywords{
Eigenvalues,
Estimation of eigenvalues, upper and lower bounds,
Variational methods for eigenvalues of operators
}

\begin{abstract}
In this article we further develop
a perturbation approach to the Rayleigh--Ritz approximations
from our earlier work. We both sharpen the estimates
and extend the
applicability of the theory to nonnegative definite operators .
The perturbation argument enables us to solve two problems in one go:
We determine which part of the spectrum of the
operator is being
approximated by the Ritz values and compute
the approximation estimates. We also present a Temple--Kato like inequality which
---unlike the original Temple--Kato inequality---
applies to any test vectors from the quadratic form domain of the operator.
\end{abstract}

\maketitle
\pagestyle{myheadings}
\thispagestyle{plain}
\markboth{On eigenvalue and eigenvector estimates for nonnegative definite operators}%
{L. Grubi\v{s}i\'{c}}

\section{Introduction}
A perturbation approach to Rayleigh--Ritz approximation was introduced by Kahan
in \cite{KahanReport}.
The main idea is to represent the eigenvalues (vectors), which we do not know
(but want to approximate), as perturbations of the Ritz values (vectors) which
we have computed.
This concept lies behind the standard subspace approximation
theory from  \cite{DavisKahan70,DavisKahWein82}. In our previous paper \cite{GruVes02}
we have shown a way to apply this concept to less regular test spaces than those
which were considered in \cite{DavisKahan70,DavisKahWein82}. In the present note
we continue this study and both improve and generalize the perturbation estimates from \cite{GruVes02}.

Let us introduce some preliminary notation. Let $h$ be a positive definite
symmetric form in a possibly infinite dimensional Hilbert space $\H$. The form $h$ generates
the positive definite operator $H$ such that $h(u,v)=(H^{1/2}u, H^{1/2}v)$.
The test space for the Rayleigh--Ritz method will be $\ra(X)$, where $X:\C^n\to\H$ is
an isometry such that $\ra(X)\subset\d(H^{1/2})$. Set $P=XX^*$, $P_\perp=\I-XX^*$ and define:
\begin{itemize}
\item the \textit{block diagonal part of} $h$ as the positive definite form\\$h'(u,v)=h(Pu, P v)+h(P_\perp u, P_\perp v)$
\item the \textit{block diagonal part of} $H$ as the operator $H'$ such that\\
$h'(u,v)=(H^{'1/2}u, H^{'1/2} v)$
\item the \textit{Rayleigh quotient} as the matrix $\Xi=(H^{1/2}X)^*H^{1/2}X\in\C^{n\times n}$.
\end{itemize}
The standard theory of \cite{DavisKahan70,DavisKahWein82} uses
\begin{equation}\label{uvod:kahan}
\max_x |(x, Hx-H^{'}x)|=\|R\|<\infty,\qquad R=HX-X\Xi=HX-H'X
\end{equation}
to obtain spectral estimates.
The operator $R$ is called the \textit{residual} of the test subspace $\ra(X)$.

It has already been demonstrated---in \cite{GruVes02}---that Kahan's concept
can yield nontrivial estimates even when $H-H'$ is not a bona fide operator,
that is to say when $\|R\|=\infty$.
We now continue the study from \cite{GruVes02} and both sharpen the estimates
and extend the applicability of the theory to
nonnegative $H$. Our results are generalizations of the known estimates for
finite matrices \cite{DrmHAri97,MathVeselic98}. A familiarity with the paper \cite{GruVes02}
is not a prerequisite for this work.

As a start we review some geometrical results from \cite{DrmHAri97}. Let for the moment
$\H$ be finite dimensional.
For the forms $h$ and $h'$ we have
\begin{align}\label{uvod:drm}
\max_x\frac{|(H^{1/2}x, H^{1/2}x)-(H^{'1/2}x, H^{'1/2}x)|}{(x, H'x)}&=\sin\Theta(H^{1/2}X,H^{-1/2}X),\\
\max_x\frac{|(H^{1/2}x, H^{1/2}x)-(H^{'1/2}x, H^{'1/2}x)|}{(x, Hx)}&=
\frac{\sin\Theta(H^{1/2}X,H^{-1/2}X)}{1-\sin\Theta(H^{1/2}X,H^{-1/2}X)},
\label{uvod:drmhari}
\end{align}
where $\sinbf\Theta(H^{1/2}X,H^{-1/2}X)$ is the sine of the maximal canonical angle
between the subspaces $H^{1/2}X$ and $H^{-1/2}X$.
We will
slightly stretch the terminology and (colloquially) call (\ref{uvod:drm}) and (\ref{uvod:drmhari})
the \textit{energy-scaled residual measures}.


Eigenvalue estimates obtained from (\ref{uvod:kahan}) are of the ``absolute'' type, i.e.
\begin{equation}\label{temple_kato_1}
|\lambda-\mu|\leq\|R\|,
\end{equation}
whereas the estimates obtained from (\ref{uvod:drm})--(\ref{uvod:drmhari})
are of the ``relative'' type
\begin{equation}\label{our_1}
\frac{|\lambda-\mu|}{\mu}\leq\sin\Theta,\qquad
\frac{|\lambda-\mu|}{\lambda}\leq \frac{\sin \Theta}{1-\sin\Theta}.
\end{equation}
We identify the following building blocks in (\ref{our_1}):
\begin{itemize}
\item $H$ and $H'$ are considered as symmetric forms \\
$h(u, v)=(H^{1/2}u, H^{1/2}v)$ and $h'(u, v)=(H^{'1/2}u, H^{'1/2}v)$
\item monotonicity of the spectrum implies the estimates
\end{itemize}
In \cite{GruVes02} the perturbation estimate (\ref{uvod:drmhari})
was shown to hold for a positive definite
operator in an infinite dimensional Hilbert space. We now prove the sharper
estimate (\ref{uvod:drm}) for a nonnegative definite operator in a Hilbert space and
give an alternative proof of (\ref{uvod:drmhari}) as a spinoff .
We also generalize ´some further results which were derived from (\ref{uvod:drmhari})
in Reference \cite{GruVes02}.

The restriction $\|R\|<\infty$, necessary for (\ref{uvod:kahan}) to give
useful information in the unbounded operator setting, incurs $\ra(X)\subset\d(H)$.
For (\ref{uvod:drm}) and (\ref{uvod:drmhari}) to be applicable we
only need to assume
$$\sin\Theta(H^{1/2}X,H^{-1/2}X)<1.$$ This ``residual measure'' will give nontrivial information
even when
$\ra(X)\subset\d(H^{1/2})$ is such that $\ra(X)\not\subset\d(H)$,
see \cite{GruVes02} and Section \ref{section_problem} of this note.

Notably, both approaches to measure the ``residual" share the property:
\begin{itemize}
\item $\sinbf\Theta(H^{1/2}X,H^{-1/2}X)=0$ if and only
if\\ $\ra(X)$ is an invariant subspace of $H$
\item $\|R\|=0$
if and only if\\ $\ra(X)$ is an invariant subspace of $H$
\end{itemize}

An important feature of our theory is that it gives an abstract framework
for a consideration of both eigenvalue and eigenvector estimates.

To get a better feeling for the estimate (\ref{our_1}) consider a simple example. Set
\begin{equation}\label{primjer_1}
H_\eta=\begin{bmatrix} \frac{1}{100} & -\frac{1}{100} \cr -\frac{1}{100} & 1 + \eta^2  \end{bmatrix}=
\begin{bmatrix} 1 & 0 \cr-1 & 1   \end{bmatrix}
\begin{bmatrix} \frac{1}{100} & 0\cr 0 &  \eta^2  \end{bmatrix}
\begin{bmatrix} 1 & -1 \cr 0 & 1   \end{bmatrix}
\end{equation}
and $e=\begin{bmatrix}1&0\end{bmatrix}^*$.
We will analyze
an approximation of the first eigenvalue of the matrix $H_\eta$ by the
Ritz value $\mu_e=e^*H_\eta e=10^{-2}$ for $\eta$ large.

As a starting point for developing a practical procedure to compute the
estimates (\ref{our_1}) we use the formula
\begin{equation}\label{intro:formula}
\sin^2\Theta(H^{1/2}X,H^{-1/2}X)=\max_{x\in\ra(X)}\frac{(x, H^{-1}x)- (x, H^{'-1}x)}{(x, H^{-1}x)},
\end{equation}
which is
implicit in \cite[Section 4.]{GruVes02}.
Since
$$
H^{'-1}=\begin{bmatrix}100&0\\0&\frac{1}{1+\eta^2}\end{bmatrix},\qquad
H^{-1}=\begin{bmatrix} 100 + {\eta }^{-2} & {\eta }^{-2} \cr {\eta }^{-2} & {\eta }^{-2}\end{bmatrix}
$$
we compute, with a help of (\ref{intro:formula}),
\begin{align*}
\lambda_1(H_\eta)&=\frac{1 + 50\,{\eta }^2 - {\sqrt{1 + 2500\,{\eta }^4}}}{100}\\
\lambda_2(H_\eta)&=\frac{1 + 50\,{\eta }^2 + {\sqrt{1 + 2500\,{\eta }^4}}}{100}\\
\sin\Theta(H^{1/2}e,H^{-1/2}e)&=\frac{1}{\sqrt{100\eta^2+1}}.
\end{align*}
As a comparison we will use an estimate which can be obtained
from the Temple--Kato inequality from \cite{ReedSimonSvi},
see (\ref{residual-temple-kato}) below.
The obtained lower bounds for $\lambda_1(H_\eta)$ are displayed
in Table \ref{f_femmod2}.
\begin{table}[htb]
\begin{center}
\begin{tabular}{|c|c|c|}\hline
  \begin{minipage}[c]{1cm}\begin{center}$~$\\[0.3em]$\eta$ \end{center}
  \end{minipage}&\begin{minipage}[c]{5cm}\begin{center}
  $~$\\[0.3em]$\text{Temple--Kato}$
  \end{center}\end{minipage}&
  \begin{minipage}[c]{5cm}\begin{center}
  $~$\\[0.3em]$
  \sin\Theta
  $
  \end{center}\end{minipage}\\[1em]\hline
  &&\\
  1  &  9.998000499860e-7         &  0.009004962810   \\
  2  &  0.007500015625   &  0.009500623831   \\
  3  &  0.008888890261   &  0.009666851698   \\
  4  &  0.009375000244   &  0.009750078088   \\
  5  &  0.009600000064   &  0.009800039988   \\
 &&\\\hline
\end{tabular}\\[1em]

\caption{Lower estimates for $\lambda_1(H_\eta)$ which can be obtained from
the Ritz value $\mu_e=10^{-2}$ with a use of
Temple--Kato estimate and with a use of $\sinbf\Theta$ approach.}\label{f_femmod2}
\end{center}
\end{table}

We can observe in Table \ref{f_femmod2}
the same behavior which was showed on an infinite dimensional model problem from \cite{GruVes02}.
Namely, the estimate $$
(1-\sin\Theta)\mu_e\leq\lambda_1,
$$
which is linear in $\sinbf\Theta$, outperforms the estimate
\begin{equation}\label{residual-temple-kato}
\mu_e-\frac{\|H_\eta e-H_\eta' e\|^2}{\lambda_2-\mu}\leq\lambda_1,
\end{equation}
which is quadratic
 in $\|H_\eta e-H_\eta' e\|$.

As an infinite dimensional analogue of (\ref{primjer_1}) we consider the following
operator.
Let $\chi_{[1,2]}$ be the characteristic function of the interval $[1,2]\subset\R$.
We consider $\mH_\eta$ which is defined by
\begin{equation}\label{e_toplina_u}
(\mH_\eta^{1/2} u,\mH_\eta^{1/2} v)=\int_{0}^2 (1+\eta^2\chi_{[1,2]})u'v'~dx.
\end{equation}
and we choose
\begin{equation}\label{drugo:e_testfunkcija_u}
u_1(x)=\begin{cases}\sqrt{2}\sin(\pi x),&0\leq x\leq 1\\
0,& 1\leq x\end{cases}
\end{equation}
as a test function. Now $u_1\in\d(\mH^{1/2}_\eta)$ but $u_1\not\in\d(\mH_\eta)$ so
neither of Temple--Kato estimates (for eigenvectors or eigenvalues) does apply since $\|\mH_\eta u_1-\mu u_1\|=\infty$.

Improved eigenvalue and eigenvector approximation estimates can be summed
up in the following procedure\footnote{Here we have
assumed we are approximating the lower end of the spectrum. Analogous procedures
can be formulated for other contiguous spectral intervals.}:
\begin{itemize}
\item Let $\mH$ be positive definite and let  $P$ be an orthogonal projection such
that $\ra(P)\subset\d(\mH^{1/2})$ and $n=\textsf{dim}~\ra(P)<\infty$.
\item If $\sinbf\Theta<1$ (as defined by (\ref{intro:formula}))
then there exist $n$-eigenvalues of the operator $\mH$ which
are approximated by the $n$ Ritz values from the subspace $\ra(P)$ in the
sense of (\ref{our_1}).
\item If $\frac{\sin\Theta}{1-\sin\Theta}<\frac{\lambda_{n+1}-\mu_n}{\lambda_{n+1}+\mu_n}$
then the Ritz values from the subspace $\ra(P)$ approximate first $n$ eigenvalues of $\mH$ (counting
the eigenvalues according to their multiplicities) and we have an eigenvector estimate.
(Analogous estimates hold for any other contiguous spectral interval.)
\end{itemize}

\section{The notation and preliminaries}
The environment in this article will be a Hilbert space $\H$, with the scalar product
$(\cdot,\cdot)$. The scalar product is antilinear in the first variable and linear in the
second. We start with a closed symmetric form $h(\cdot,\cdot)$ which is
additionally assumed to be \textit{nonnegative}\index{form!nonnegative}
\begin{equation}\label{druga:nonnegative}
h[u]=h(u,u)\geq 0,\qquad u\in\q(h).
\end{equation}
In the sequel when we say the nonnegative form $h$, we shall always mean the closed symmetric form $h$
which satisfies (\ref{druga:nonnegative}). The form $h$ shall be called \textit{positive definite}
\index{form!positive definite}
when it is closed, symmetric and there exists $m_h>0$ such that
$$
h[u]=h(u,u)\geq m_h\|u\|^2 ,\qquad u\in\q(h).
$$
There is also an equivalent operator version of these definitions.
The selfadjoint operator $\mH$ is called
\textit{nonnegative}\index{operator!nonnegative} if
$$
(u,\mH u)\geq 0,\qquad u\in\d(\mH).
$$
Subsequently, $\mH$ is called \textit{positive definite}
\index{operator!positive definite}
if there exists $m_\mH>0$ such that
$$
(u,\mH u)\geq m_\mH\|u\|^2 ,\qquad u\in\d(\mH).
$$
In this chapter we assume $\overline{\q}^{~_\H}=\H$, but later we shall also allow
$\overline{\q}^{~_\H}$ to be any nontrivial subspace of $\H$. For nonnegative
selfadjoint operators one defines, with the help of the spectral theorem,
the usual functional calculus. We write the spectral decomposition of the
selfadjoint operator $\mH$ as
$$
\mH=\int\lambda~\text{d}E_\mH(\lambda),
$$
where $E_\mH(\lambda)$ is the right continuous \textit{spectral family}
associated to the operator $\mH$. When there can be no confusion we simply write
$E(\lambda)$.

The representation theorem for
nonnegative forms \cite[pp. 331]{Kato76} implies  that there exists a selfadjoint
operator $\mH$ such that $\d(\mH^{1/2})=\q(h)$ and
$$
h(u,v)=(\mH^{1/2}u,\mH^{1/2}v),\qquad u,v\in\q(h)
$$
Following \cite{Faris75} we call $\d(\mH)$ the \emph{operator domain}\index{domain!operator domain}
of $\mH$ and $\q(\mH)=\d(\mH^{1/2})$ the \emph{quadratic form
domain}\index{operator!quadratic form domain} of $\mH$. We write $\d$ and $\q$ when
there can be no confusion. With the help of the spectral theorem we see that
\begin{align*}
\d(\mH)&=\{u\in\H:\|\mH u\|^2=\int\lambda^2~d(E(\lambda)u,u)<\infty\},\\
\q(\mH)&=\{u\in\H:h[u]=\|\mH^{1/2} u\|^2=\int\lambda~d(E(\lambda)u,u)<\infty\}.
\end{align*}

In general, when dealing with the forms in a Hilbert space we
shall follow the terminology of Kato, cf. \cite{Kato76}. In one
point we will depart from the conventions in \cite{Kato76}.
A nonnegative form $$h(u,v)=(\mH^{1/2}u, \mH^{1/2}v)$$
will be called \emph{nonnegative definite}\index{form!nonnegative definite}
when  $\lambda_e(\mH)>0$. Analogously, a nonnegative operator $\mH$ such that
$\lambda_e(\mH)>0$ will be also called
\emph{nonnegative definite}\index{operator!nonnegative definite}.
We will often say nonnegative, meaning the nonnegative definite.
Now, we give
definitions of some terms that will frequently be used, cf. \cite{Faris75,Kato76}.

\begin{definition}
A bounded operator $A:\H\to\mathcal{U}$ is called \textit{degenerate}\index{operator!degenerate}
if $\ra(A)$ is finite dimensional.
\end{definition}

\begin{definition}
Let $\mH$ and $\mA$ be nonnegative operators. We define the \textit{order relation}
\index{operator!the order relation}
$\leq$ between the nonnegative operators by saying that
$$
\mA\leq\mH
$$
if $\q(\mH)\subset\q(\mA)$ and
$$
\|\mA^{1/2} u\|\leq\|\mH^{1/2} u\|,\qquad u\in\q(\mH),
$$
or equivalently if
$$a[u]\leq h[u], \qquad u\in\q(h),$$
when $a$ and $h$ are nonnegative forms defined by the operators $\mA$ and
$\mH$ and $\mA\leq\mH$.
\end{definition}

A main principle we shall use to develop the perturbation theory will be
the \textit{monotonicity of the spectrum}\index{spectrum!monotonicity of the spectrum}
with regard to the order relation
between nonnegative operators. This principle can be expressed in many ways.
The relevant results, which are scattered over the monographs \cite{Faris75,Kato76}, are summed
up in the following theorem, see also \cite[Corollary A.1]{Levendorskii}.

\begin{theorem}\label{prvo:t_monotonost}
Let $\mA=\int\lambda~\text{d}E_\mA(\lambda)$ and $\mH=\int\lambda~\text{d}E_\mH(\lambda)$
be nonnegative operators in $\H$ and let $\mA\leq\mH$. By
$0\leq\mu_1\leq\mu_2\leq\cdots<\lambda_e(\mA)$ and
$0\leq\lambda_1\leq\lambda_2\leq\cdots<\lambda_e(\mH)$
denote the discrete eigenvalues of $\mA$ and $\mH$, then
\begin{enumerate}
\item $\lambda_e(\mA)\leq\lambda_e(\mH)$
\item ${\sf dim} ~E_\mH(\gamma)\leq {\sf dim} ~E_\mA(\gamma)$, for every $\gamma\in\R$
\item $\mu_k\leq\lambda_k,\qquad k=1,2,\cdots$ .
\end{enumerate}
\end{theorem}

We close this introductory section with the well known theorem about
the perturbation of the \textit{essential spectrum}\index{spectrum!essential spectrum}.

\begin{theorem}\label{prvo:t_weyl}
Let $\mH$ and $\mA$ be positive definite operators. If the operator
$$
\mH^{-1}-\mA^{-1}
$$
is compact then $\sigma_{ess}(\mH)=\sigma_{ess}(\mA)$.
\end{theorem}

\section{The generalized inverse and angle between the subspaces}
\label{prvo:s_kutevi}

There are many ways to express that $u\in\q(h)$ is an eigenvector
of the operator $\mH$. We will give a geometric characterization
of this property. Assume that $\|u\|=1$ and $\mu=h[u]$. An elementary
trigonometric argument yields
\begin{equation}\label{prvo:eq:elementary}
\|\mH^{1/2}u-\mu \mH^{-1/2}u\|=0 \Leftrightarrow \sin\Theta(\mH^{1/2}u, \mH^{-1/2}u)=0.
\end{equation}
 (\ref{prvo:eq:elementary}) implies that
$u$ is an eigenvector of $\mH$ if and only if $\sinbf\Theta(\mH^{1/2}u, \mH^{-1/2}u)=0$.
The ability to assess the size of
$\sinbf\Theta(\mH^{1/2}u, \mH^{-1/2}u)$ will
be central to the analysis of the Rayleigh--Ritz method in
this paper.

In this section we give the background information on the angles between
two finite dimensional subspaces of a Hilbert space as given in \cite{DavisKahan70,Kato76,wed-82}.
Basic results on generalized inverses of (unbounded) operators defined between two
Hilbert spaces will be presented as well. These results will be applied
to the problem of computing $\sinbf\Theta(\mH^{1/2}\mathcal{X},\mH^{-1/2}\mathcal{X})$
for the given positive definite $\mH$
and some finite dimensional $\x\subset\q(\mH)$.

Closed subspaces of the Hilbert space $\H$ can be represented as images
of the corresponding orthogonal projections. We shall freely speak about
the dimension of the projection $P$ meaning the dimension
of the range of the projection $P$. In the case in which $P$ is finite dimensional, we
have another representation for the subspace $\ra(P)$.
For a given $n$-dimensional subspace $\ra(P)\subset\q$ there exists an isometry $X:\C^n\to\H$
such that $\ra(P)=\ra(X)$, where $P=XX^*$. Therefore, $\ra(X)$ is an alternative representation of
the $n$-dimensional subspace $\ra(P)$. The isometry $X$ will be called the
basis of the subspace $\ra(P)$. We shall freely use both representation of
the finite dimensional subspace. $P_X=XX^*$ will generically denote the
orthogonal projection on the space $\ra(X)$ (for some isometry $X:\C^n\to\H$).

Let $\ra(P)$ and $\ra(Q)$ be two finite\index{basis of a subspace}
dimensional subspaces of the Hilbert space $\H$. The function $\angle$ that measures
the separation of the pair of subspaces $\ra(P)$ and $\ra(Q)$ will be called
an \textit{angle function}\index{angle function} if it
satisfies the following  properties
\begin{enumerate}
    \item $\angle(P,Q)\geq 0$ and \\$\angle(P,Q)=0$ if and
    only if $\ra(P)\subset \ra(Q)$ or $\ra(Q)\subset\ra(P)$.
    \item $\angle(P,Q)=\angle(Q,P)$
    \item $\angle(P,Q)\leq\angle(P,R)+\angle(R,Q)$ if \\${\sf dim}(\ra(P))\leq
    {\sf dim}(\ra(R))\leq{\sf dim}(\ra(Q))$ or\\
    ${\sf dim}(\ra(P))\geq
    {\sf dim}(\ra(R))\geq{\sf dim}(\ra(Q))$
    \item $\angle(UP,UQ)=\angle(P,Q)$, for any unitary $U$.
\end{enumerate}

In what follows we will use the following angle functions, see \cite{wed-82},
\begin{align}\label{eq:sinTheta}
\Theta(P,Q)&=\arcsin \max\{\|P(\I-Q)\|,\|Q(\I-P)\|\}\\
\Theta_p(P,Q)&=\arcsin \min\{\|P(\I-Q)\|,\|Q(\I-P)\|\}\label{eq:principalsinT}
\end{align}
The function $\Theta(P,Q)$ from (\ref{eq:sinTheta})
will be called the \textit{maximal canonical angle}\index{angle!maximal canonical angle}
between the subspaces $P$ and $Q$. The function $\Theta_p(P,Q)$ from (\ref{eq:principalsinT})
will be called the \textit{maximal principal angle}\index{angle!maximal principal angle} between the subspaces $P$
and $Q$.

The following lemma, which is a consequence of \cite[Theorem I-6.34]{Kato76},
gives an insight in the behavior of the canonical and
the principal angles which were defined by (\ref{eq:sinTheta}) and (\ref{eq:principalsinT}).

\begin{lemma}\label{prvo:hoce_sef}
Let $P$ and $Q$ be two orthogonal projections such that ${\sf dim}(\ra(P))\leq{\sf dim}(\ra(Q))$
and let
$$
\|P(\I-Q)\|<1
$$
then we have the following alternative. Either
\begin{enumerate}
    \item ${\sf dim}(\ra(P))={\sf dim}(\ra(Q))$ and
    $$
    \sin\Theta(P,Q)=\sin\Theta_p(P,Q)=\|P-Q\|<1,\qquad\text{or}
    $$
    \item ${\sf dim}(\ra(P))<{\sf dim}(\ra(Q))$ and
    $$
    \sin\Theta_p(P,Q)=\|P(\I-Q)\|<1.
    $$
\end{enumerate}

\end{lemma}

For most of our needs, Lemma \ref{prvo:hoce_sef} describes the relation between the
finite dimensional subspaces $\ra(P)$ and $\ra(Q)$ in
sufficient detail. However, sometimes
it will be necessary to analyze the structure of
the finite dimensional projections $P_V=VV^*$ and $P_U=UU^*$
in further detail. To this end we define
the \textit{canonical angles}\index{angle!canonical angles} $\theta_1,\ldots ,\theta_n$
between the spaces $\ra(U)$ and $\ra(V)$
 as
\begin{equation}\label{prvo2:e_defkutevi}
\theta_i=\arccos\sigma_i,\quad i=1,\ldots,n,
\end{equation}
where $\sigma_1,\ldots,\sigma_n$ are the singular values of the matrix
$$
V^*U\in\C^{m\times n}.
$$
We have assumed that $m=\textsf{dim}~\ra(V)$, $n=\textsf{dim}~\ra(U)$ and $m\leq n$.
The canonical angles are related to the angle function (\ref{eq:sinTheta}) through the
formula, see \cite{wed-82},
$$
\sin\Theta(P_V,P_U)=\max_i\sin\theta_i.
$$
We also define the \textit{acute principal angles} $\theta^p_1\leq\theta^p_2\leq\cdots\leq\theta^p_k$, where
$k\leq n$,
\index{angle!acute principal angles} as those canonical angles $\theta_i$ which
satisfy the condition $0<\theta_i<\pi/2$. Subsequently, we obtain a connection to the
angle function (\ref{eq:principalsinT}) through the formula
$$
\sin\Theta_p(P_V,P_U)=\max_i\sin\theta^p_i.
$$

In dealing with the projections and degenerate operators it is useful
to have a notion of the generalized inverse. We will use the definition of the generalized
inverse of a closed densely defined operator in $\H$ from \cite{Nashed76}, see
also \cite[Chapter IV.5]{Kato76}.

\begin{definition}
Let $\mathbf{T}:\H\to\mathcal{U}$ be a closed operator such
that $\overline{\d(\mathbf{T})}=\H$. The operator $\mT^\dagger:\mathcal{U}\to\H$ is
defined by
\begin{align*}
\d(\mathbf{T}^\dagger)&=\ra(\mathbf{T})\oplus\ra(\mathbf{T})^\perp\\
\mT^\dagger u&=(\left.\mT_{~}\right|_{\je(\mT)^\perp})^{-1}P_{\ra(\mT)}u,\qquad
u\in\d(\mT^\dagger)
\end{align*}
and it is called the \textit{Moore-Penrose generalized inverse}
\index{generalized inverse!More-Penrose generalized inverse} of $T$.
\end{definition}

The properties of the generalized inverse\footnote{The generalized
inverses can also be defined in more general settings. Their properties are also analyzed
in \cite{Nashed76}.} are analyzed
in the monograph \cite{Nashed76}. In particular we use the following characterization.
\begin{theorem}[see {\cite[Theorem I.5.7]{Nashed76}}]\label{prvo:theorem:1}
Let $\mathbf{T}:\H\to\mathcal{U}$ be the closed operator and let
$\overline{\d(\mathbf{T})}=\H$, then $\mathbf{T}^\dagger$ is the unique closed operator such that
\begin{align*}
\mathbf{T}^\dagger\mathbf{T}\mathbf{T}^\dagger&=\mathbf{T}^\dagger,
\qquad \text{on } \d(\mathbf{T}^\dagger)\\
\mathbf{T}\mathbf{T}^\dagger&=\left.P_{\ra(\mathbf{T})}\right|_{\d(\mathbf{T}^\dagger)}\\
\mathbf{T}^\dagger\mathbf{T}&=\left.P_{\je(\mathbf{T})^\perp}\right|_{\d(\mathbf{T})}
\end{align*}
where $P_{\mathcal{M}}$ is the orthogonal projection on $\mathcal{M}$.
The operator $\mT^\dagger$ is bounded if and only if $\mT$ has a closed range.
\end{theorem}

The nonnegative operator $\mH^\dagger$ has
the spectral decomposition
$$
\mH^{\dagger}=\int\frac{1}{\lambda}~dE(\lambda),\qquad \d(\mH^\dagger)=\{u\in\H:\int\frac{1}{\lambda^2}
~d(E(\lambda)u,u)<\infty\},
$$
and the functional calculus implies
$$
\mH^{\dagger 1/2}=\mH^{1/2\dagger}.
$$

Theorem \ref{prvo:theorem:1} shows a relation between
the Moore-Penrose generalized inverses and orthogonal projections in
a Hilbert space. This is precisely the reason why the generalized inverses
will be useful in our study.

A bounded operator $W:\H\to\mathcal{U}$ is called \textit{partially isometric}
\index{operator!partialy isometric}
if there exists a closed subspace $\mathcal{M}\subset\H$ such that
$$
\|Wu\|=\|P_{\mathcal{M}}u\|,\qquad u\in\H.
$$
This is equivalent to
$$
W^*W=P_{\mathcal{M}}.
$$
The set $\mathcal{M}=\ra(W^*)\subset\H$ is called the \textit{initial set}
 of the partial isometry $W$ and
$\ra(W)\subset\mathcal{U}$ is called the \textit{final set}. Since $\je(W^*)\oplus\ra(W)$ we
see
$$
WW^*=P_{\ra(W)},
$$
so $W^*$ is also the partial isometry with the initial set $\ra(W)$.
We shall also use the notation
\begin{align*}
W^*W&=P_{W^*}, &WW^*=P_W.
\end{align*}
It is obvious that
$$
W^*=W^\dagger
$$
and we have the following lemma.

\begin{lemma}\label{prvo:lema:parcijalneizometrije}
A bounded operator $W:\H\to\mathcal{U}$ is partially isometric if and only if
$$
WW^*W=W.
$$
\end{lemma}
For the proof see \cite{Kato76}.
\begin{lemma}\label{prvo1:l_sinuskuta}
Let $V$ and $W$ be two partial isometries then
$$
\|P_VP_W\|=\|VP_W\|=\|V^*W\|.
$$
\end{lemma}
\begin{proof}
Using Lemma \ref{prvo:lema:parcijalneizometrije} we compute
\begin{align*}
\|P_VP_W\|^2&=\text{spr}(P_WP_VP_W)=\text{spr}(WW^*VV^*WW^*)\\
&=\text{spr}(W^*VV^*WW^*W)=\text{spr}(W^*VV^*W)=\|V^*W\|.
\end{align*}
In this computation we have used the identity
$$
\text{spr}(ABC)=\text{spr}(CAB),
$$
which holds for bounded operators $A, B, C$.
\qquad\end{proof}


\section{Geometrical properties of the Ritz value perturbation}\label{prvo:s:perturnbation}

In this section we will present a perturbation approach to
the Rayleigh--Ritz approximation of the spectrum of a positive definite operator.
The nonnegative definite case is technically more complex and warrants a separate section.
Although this chapter is devoted to the positive definite case,
some of the statements and definitions will be given in full generality in
which they will be later used in the text.

Let $0\leq h$ be a nonnegative form and let
$\ra(X)\subset\q(h)$ be the $n$-dimensional test space. The matrix
$$
\Xi_{\mH, X}=(\mH^{1/2}X)^*\mH^{1/2}X\in\C^{n\times n}
$$
will be called the \textit{Rayleigh quotient}\index{Rayleigh quotient} associated to the basis $X$.
When there can be no confusion, we shall denote
the Rayleigh quotient by $\Xi$ and drop the indices.
The eigenvalues of the matrix $\Xi$ will be numbered in the ascending order
\begin{equation}\label{prvo:e_uredjRitz}
\mu_1\leq\mu_2\leq\cdots\leq\mu_n.
\end{equation}
We call the numbers $\mu_i$ the Ritz values of the operator $\mH$ (form $h$) from
the subspace $\ra(X)$.
This definition is correct since the eigenvalues of the
matrix $\Xi$ do not depend on the choice of
the basis $X$. In the rest of this chapter we will use $P=XX^*$
to denote the projection onto the range of the isometry
$X:\C^n\to\H$.

For the given $h$ and $\ra(X)\subset\q(h)$, $P=XX^*$, we define the symmetric forms $\delta h$ and $h'$
using the formulae
\begin{align}\label{prvo:e_def_dh}
\delta h(u,v)&=h(Pu, (\I-P)v) + h((\I-P)u, Pv),\qquad u,v\in\q(h)\\
h'(u,v)&=h(Pu, Pv)+h((\I-P)u,(\I-P)v),\qquad u,v\in\q(h).\label{prvo:e_def_h'}
\end{align}
Obviously, (\ref{prvo:e_def_dh}) and (\ref{prvo:e_def_h'}) imply
\begin{equation}\label{prvo1:e_generickiH}
h'(u,v)=h(u,v)-\delta h(u,v),\qquad u,v\in\q(h).
\end{equation}
Before we can proceed we need the following definition.

\begin{definition}\label{prvo:d_komutira}
If $\mH$ is a selfadjoint operator and $P$ a projection, to say that
$P$ commutes with $\mH$ means that $u\in\d(\mH)$ implies $Pu\in\d(\mH)$
and
$$
\mH Pu=P\mH u,\qquad u\in\d(\mH).
$$
\end{definition}

In what follows we will describe the properties of the symmetric form $h'$ and
of the operator $\mH'$ it generates.

\begin{lemma}\label{prvo1:l_reducira}\label{prvo:l_lemma-zatvorena}
Let the nonnegative definite form $h$ and the subspace $\ra(X)\subset\q$ be given.
Let $\mH$ be the nonnegative definite operator defined by the form $h$.
The form $h'$ from (\ref{prvo:e_def_h'}) is closed and positive and it defines the selfadjoint operator $\mH'$.
Furthermore, $\mH'$ is positive definite if $\mH$ is positive definite,
$\sigma_{ess}(\mH)=\sigma_{ess}(\mH')$ and
\begin{equation}\label{prvo:t_matricnoRed}
\mH'X=X \Xi,
\end{equation}
for $\Xi=(\mH^{1/2}X)^*\mH^{1/2}X\in\C^{n\times n}$.
\end{lemma}
\begin{proof}
The operators $\mH^{1/2}P$ and $\mH^{1/2}(\I-P)$ are closed and so is the form
$$
h'(u,v)=h(Pu, Pv)+h((\I-P)u,(\I-P)v).
$$
It is obviously nonnegative, so it defines a nonnegative selfadjoint operator $\mH'$.
We will now show that the subspace $\ra(X)$ reduces $\mH'$. Indeed,
for $y\in\q$, $x\in\C^n$ we have
\begin{eqnarray*}
h'(y,X x)&=&(\mH^{1/2} y, \mH^{1/2}X x) -(\mH^{1/2}(\I-P)y, \mH^{1/2}X x) \\
&=&(\mH^{1/2}XX^*y, \mH^{1/2}X x) \\
&=&(\Xi X^*y, x).
\end{eqnarray*}
This is equivalent to
$$
(\mH^{'1/2}y, \mH^{'1/2} X x) = (y, X\Xi x) ,\quad y\in\q,\;x\in\C^n.
$$
It implies $\ra(X)\subset\d(\mH')$ and
$$
(y, \mH' X x- X\Xi x)=0,
$$
for all $y\in\H$, $x\in\C^n$. Hence,
\begin{equation}\label{prvo1:e_jednreducira}
\mH'X=X\Xi
\end{equation}
which is equivalent to the statement that $P$ commutes with $\mH'$
(see Definition \ref{prvo:d_komutira}).
We now prove that $\sigma_{ess}(\mH)=\sigma_{ess}(\mH')$.
Assume $h$ is a positive definite form, then
 $h'$ from (\ref{prvo:e_def_h'}) is positive definite, too. From (\ref{prvo1:e_generickiH}) we obtain
$$
\delta h(\mH^{-1}u,\mH^{'-1}v)=(\mH^{'-1}u-\mH^{-1}u,v),\qquad u,v\in\H.
$$
On the other hand
$$
\delta h(\mH^{-1}u,\mH^{'-1}v)=(\mH^{1/2}P\mH^{-1}u,\mH^{1/2}P_\perp\mH^{'-1}v)+
(\mH^{1/2}P_\perp\mH^{-1}u,\mH^{1/2}P\mH^{'-1}v)
$$
defines a compact operator. Theorem \ref{prvo:t_weyl} implies
$\sigma_{ess}(\mH)=\sigma_{ess}(\mH')$ and the statement of
the theorem is proved for a positive definite $h$. In the general case, take $\alpha>0$.
The form $\widetilde{h}(u, v)=h(u, v)+\alpha(u, v)$ is positive definite. Furthermore,
 we establish
\begin{align*}
\widetilde{h}'(u,v)&=\alpha (u, v)+h'(u,v)\\
\delta\widetilde{h}(u, v)&=\delta h(u, v),
\end{align*}
so $\sigma_{ess}(\widetilde{\mH})=\sigma_{ess}(\widetilde{\mH}')$. The conclusion
$\sigma_{ess}(\mH)=\sigma_{ess}(\mH')$ follows by the
spectral mapping theorem. \qquad\end{proof}

\begin{corollary}\label{prvo:corr:333}
Let the nonnegative definite form $h$ and the subspace $\ra(X)\subset\q$ be given.
The projections $P=XX^*$ and $P_{\ra(\mH')}$ commute and $\je(\mH')\subset\je(\mH)$.
\end{corollary}

\begin{remark}{\rm
For positive definite $h$ Lemma \ref{prvo1:l_reducira} describes
the operator $\mH'$ in sufficient detail. For a general
nonnegative $h$ the operator $\mH'$ has somewhat more complex structure.
Further properties of the operator $\mH'$, constructed in the case in which
$h$ is a general nonnegative form, will be discussed in Section \ref{prvo1:s_nonnegativecase}.}
\end{remark}

We now concentrate on the positive definite case.

\begin{theorem}\label{prvo1:l_lema13}
Let the subspace $\ra(X)\subset\q$ be given and let $h$ be positive definite.
Assume $\sinbf\Theta:=\sinbf\Theta(\mH^{1/2} X,\mH^{-1/2}X)<1$, then
\begin{align}\label{prvo:e_prvizakj}
(1-\sin\Theta)h'[u]
&\leq h[u]\leq(1+\sin\Theta)h'[u],\qquad u\in\q(h)\\
(1-\frac{\sin\Theta}{1-\sin\Theta})h[u]
&\leq h'[u]\leq(1+\frac{\sin\Theta}{1-\sin\Theta})
h[u],\qquad u\in\q(h).\label{prvo:e_drugizakj}
\end{align}
\end{theorem}
\begin{proof}
The product $\mH^{1/2}\mH'^{-1/2}$ is well
defined since $\q=\d(\mH^{1/2})=\d(\mH'^{1/2})$. This implies that the
form
$$
\delta h_s(x,y)=\delta h(\mH^{'1/2}x, \mH^{'1/2} y)
$$
defines the bounded operator $\delta H_s$.
After the substitutions $u=\mH'^{-1/2} x$, $v=\mH'^{-1/2} y$ we obtain
\begin{equation}\label{prvo:e_sinVW}
\max_{u, v\in\q(h)}\frac{|\delta h(u, v)|}{\sqrt{h'[u] h'[v]}}=\|\delta H_s\| .
\end{equation}
We now show $\|\delta H_s\|=\sinbf\Theta$.
Set\begin{align}\label{prvo1:eq_WV_1}
V&=\mH^{1/2}P\mH^{'-1/2}\\
W&=\mH^{1/2}P_\perp\mH'^{-1/2},
\end{align}
with $P_\perp=\I-P$.
Relation (\ref{prvo1:e_generickiH}) implies
\begin{align}
\nonumber \delta h(\mH^{'-1/2}u, \mH^{'-1/2}v)&=h(P_\perp \mH^{'-1/2}u,P\mH^{'-1/2}v)
+ h(P\mH^{'-1/2}u, P_\perp\mH^{'-1/2}v)\\
&=(Wu, Vv)+(Vu, Wv),\label{prvo1:eq_WV_2}
\end{align}
which can be written as
\begin{equation}\label{eq:formula_Hs}
\delta H_s=V^*W+W^*V.
\end{equation}
The equations (\ref{prvo1:eq_WV_1})--(\ref{eq:formula_Hs}) yield
\begin{align}
VW^*&=WV^*=0\\
\|\delta H_s\|&=\|W^*VV^*W+V^*WW^*V\|=\|V^*W\|.
\end{align}
As the next step we establish that $V$ and $W$ are partial isometries such that
\begin{align}\label{prvo1:eq_trivijalna}
\ra(V)&=\ra(\mH^{1/2}P)\\
\ra(W)^\perp&=\ra(\mH^{-1/2}P).\label{prvo1:eq_nijetrivijalna}
\end{align}
%

The proof will follow from Lemma \ref{prvo1:l_reducira}. It runs along the same lines in
both cases, so we will only present the proof for $W$. Take some
$u,v\in\H$, then
\begin{align*}
(W u, W v)&=(\mH^{1/2}P_{\perp}\mH'^{-1/2}u, \mH^{1/2}P_{\perp}\mH'^{-1/2}v)\\
&=h(P_{\perp}\mH'^{-1/2}u,P_{\perp}\mH'^{-1/2}v)=h'(P_{\perp}\mH'^{-1/2}u,
P_{\perp}\mH'^{-1/2}v)=(P_{\perp}u, v),
\end{align*}
so $W^*W=P_{\perp}$. This proves that $W$ is a partial isometry.

Relation (\ref{prvo1:eq_trivijalna}) is obvious, since
$$
\ra(\mH^{1/2}P\mH^{'-1/2})=\ra(\mH^{1/2}P)
$$
is guaranteed by the assumption $\ra(P)\subset\q(h)$ and
the injectivity of $\mH^{'-1/2}$.

The proof of (\ref{prvo1:eq_nijetrivijalna}) requires a bit more
work. One computes
$$
W^*\mH^{-1/2}P=\mH^{'-1/2}P_\perp\mH^{1/2}\mH^{-1/2}P=0,
$$
which implies
$$
\ra(\mH^{-1/2}P)\subset\je(W^*)=\ra(W)^\perp.
$$
On the other hand \begin{equation}\label{eq:homeomorfizam}
W^*=P_\perp A,
\end{equation} where
$A=\overline{\mH^{'-1/2}\mH^{1/2}}:\H\to\H$ is a homeomorphism
(of linear topological vector spaces), so
$$
{\sf dim}~\je(W^*)={\sf dim}~\je(P_\perp)={\sf dim}~\ra(P)={\sf dim}~\ra(\mH^{-1/2}P)
$$
and (\ref{prvo1:eq_nijetrivijalna}) is established.
The assumption $\sinbf\Theta<1$ and Lemma \ref{prvo1:l_sinuskuta} guarantee
$$
\sin\Theta=\|V^*W\|.
$$
Finally, using (\ref{prvo:e_sinVW}) we establish
$$
(1-\sin\Theta)h'[u]\leq h[u]\leq(1+\sin\Theta)h'[v],
$$
which is the statement (\ref{prvo:e_prvizakj}).

It is a well
known fact that given some $0<\lambda,\mu$ and $0<\eta<1$ the
implication
\begin{equation}\label{prvo1:e_etadef1}
\frac{|\lambda-\mu|}{\mu}\leq\eta\Rightarrow\frac{|\lambda-\mu|}{\lambda}\leq\frac{\eta}{1-\eta}
\end{equation}
holds. Since $h$ and $h'$ are positive definite forms,
the relation (\ref{prvo:e_drugizakj})
is proved.
\qquad\end{proof}

\begin{example}{\rm
Let $-\partial_{xx}$ be considered as the selfadjoint operator
with $$\d(-\partial_{xx})=\{u\in H^2[0,1] :~u(0)=u(1)=0\}.$$ The
partial integration establishes that $-\partial_{xx}$ is defined
by the positive definite form
\begin{equation}\label{prvo:e_primjerzaR}
h(u,v)=\int^1_0 \partial_x u~\partial_x v~d x,\qquad u,v\in \q(-\partial_{xx})=H^1_0[0,1].
\end{equation}
The operator $\partial_x u$, $u\in H^1_0[0,1]$ is closed, but
not selfadjoint, therefore (\ref{prvo:e_primjerzaR}) is
an alternative operator representation (factorization), to the ``square root"
representation (\ref{prvo:e_repza1_2})
of the form $h$ (the operator $-\partial_{xx}$).}
\end{example}

Take any
positive definite form $h$, then
\begin{equation}\label{prvo:e_repza1_2}
h(u,v)=(\mH^{1/2}u,\mH^{1/2}v)
\end{equation}
is only one of the possible operator representations of the form $h$. All of the preceding
results are independent of the choice of the operator representation $h(u,v)=(\mR u, \mR v)$, since
\begin{equation}\label{eq:corollary_main}
\sin\Theta=\max_{u,v\in\q}\frac{|\delta h(u,v)|}{\sqrt{h'[u]h'[v]}}
\end{equation}
and $h'$ depends only on $h$ and $\ra(P)$.

Furthermore, all of the representations of the form $h$ are in a sense equivalent.
Let $\mR:\H\to\H'$ be a closed operator such that
\begin{equation}\label{prvo:e_alternativna}
h(x,y)=\sk{\mR x,\mR y}=\sk{\mH^{1/2}x,\mH^{1/2}y}
\end{equation}
and $\q=\d(\mR)=\d(\mH^{1/2})$, then by \cite[Ch. VI.7]{Kato76}
\begin{equation}\label{prvo2:e_izometrija1}
\mR=U\mH^{1/2}, \quad \mR^*=\mH^{1/2}U^*,
\end{equation}
where $U$ is the isometry from $\H'$ onto $\ra(\mR)$.
Independence of the estimate (\ref{prvo:e_prvizakj}) from the representation
(\ref{prvo:e_alternativna}) could have also been proved by the unitary invariance of
the canonical angle and (\ref{prvo2:e_izometrija1}).
Formula (\ref{eq:corollary_main}) is an important corollary
of Theorem \ref{prvo1:l_lema13}.
In the next theorem we prove that also,
\begin{equation}\label{eq:drmac_corollary}
\frac{\sin\Theta}{1-\sin\Theta}=\max_{u,v\in\q}\frac{|\delta h(u,v)|}{\sqrt{h[u]h[v]}}
\end{equation}
holds. Equations (\ref{eq:corollary_main}) and (\ref{eq:drmac_corollary}) demonstrate
that the constants $\sinbf\Theta$ and $\frac{\sinbf\Theta}{1-\sinbf\Theta}$ in
(\ref{prvo:e_prvizakj}) and (\ref{prvo:e_drugizakj}) cannot be improved upon.

The following lemma is taken out of the
joint paper \cite{GruVes02}, cf. \cite{RitzDrm96}. We present it here without a proof.

\begin{lemma}\label{prvo2:t_drmac}
Let the form $h$ be positive definite and let
the forms $h'$ and $\delta h$ be as in (\ref{prvo1:e_generickiH}), then
\begin{equation}\label{drugo:maximum_r_h}
\max_{u,v\in\q}\frac{|\delta h(u,v)|}{\sqrt{h[u]h[v]}}=\frac{\sin\Theta}{1-\sin\Theta}
\end{equation}
holds.
Here $\sinbf\Theta=\sinbf\Theta(\mH^{1/2}X,\mH^{-1/2}X)$, where $\ra(X)\subset\q$
was the subspace used to define $h'$ and $\delta h$.
\end{lemma}


\subsection{The nonnegative definite case}
\label{prvo1:s_nonnegativecase}

In the nonnegative case we have to provide an alternative
definition for a subspace that will play the
role of $\ra(\mH^{-1/2}X)$. We have shown
$W=\mH^{1/2}P_{\perp}\mH'^{-1/2}$ to be a partial isometry
such that
$$
\mathcal{W}=\ra(\mH^{1/2}P_{\perp})^\perp=\ra(W)^\perp=\ra(\mH^{-1/2}X).
$$
The left part of the equality is also well defined in the case in which
$\mH^{1/2}$ is not invertible, so we set
$$
\mathcal{W}=\ra(\mH^{1/2}P_\perp)^\perp.
$$

The construction (\ref{prvo1:e_generickiH}) was performed with the assumption
that $h$ is nonnegative definite and ${\ra(X)}\subset\q$. Lemma \ref{prvo1:l_reducira}
says $\sigma_{ess}(\mH)=\sigma_{ess}(\mH')$ so $\mH^{'\dagger 1/2}$ is a bounded operator
and
\begin{align}\label{prvo:e_defV}
V&=\mH^{1/2}P\mH^{'\dagger 1/2},\\
W&=\mH^{1/2}P_\perp\mH^{'\dagger 1/2},\label{prvo:e_defW}
\end{align}
are everywhere defined. Corollary \ref{prvo:corr:333}
enables us to conclude that $\ra(V)=\ra(\mH^{1/2}P)$ and
$\ra(W)=\ra(\mH^{1/2}P_\perp)$, so we set
\begin{equation}\label{drugo:eq:sefubacio}
\mathcal{V}=\ra(V),\qquad \mathcal{W}=\ra(W)^\perp.
\end{equation}

Lemma \ref{prvo:l_lemma-zatvorena} states
that given a positive definite $\mH$ the constructed operator $\mH'$
must always be positive definite.
In general nonnegative situation we have only the result of Corollary \ref{prvo:corr:333}.
It establishes that $\mH'$ is a nonnegative definite operator and that
$\je(\mH')\subset\je(\mH)$.
This does not give sufficient information on the structure of $\mH'$.
Formulae like (\ref{prvo:e_prvizakj})--(\ref{prvo:e_drugizakj}) are meaningful
in the nonnegative definite case, too. They, however, invariably imply
$\je(\mH)=\je(\mH')$. We, therefore, proceed is two steps. Firstly,
we establish a general (theoretical) condition on the
subspace $\x=\ra(P)$ which guarantees
that $\je(\mH)=\je(\mH')$. As the second step we give a
practical computational formula.

The subspaces $\mathcal{W}$ and $\mathcal{V}$ need not have the same dimension,
so we will have to use the principal angle to compare them, cf. Lemma \ref{prvo:hoce_sef}.
In what follows we show that
$$
\sin\Theta_p(\mathcal{V}, \mathcal{W})
$$
takes the role of $\sinbf\Theta(\mH^{1/2}X, \mH^{-1/2}X)$ in the
nonnegative version of Theorem \ref{prvo1:l_lema13}.
In the case when $\mH^{1/2}$ is invertible (\ref{prvo1:eq_nijetrivijalna})
implies $\mathcal{V}=\ra(\mH^{1/2}X)$ and
$\mathcal{W}=\ra(\mH^{-1/2}X)$. The subspaces $\mH^{-1/2}X$ and $\mH^{1/2}X$
have the same dimension, so Corollary \ref{prvo:hoce_sef} yields
$$
\sin\Theta_p(\mathcal{V}, \mathcal{W})=\sin\Theta(\mH^{1/2}X, \mH^{-1/2}X).
$$

 We establish the properties of $V$ and $W$ and
give a characterization of the subspace $\mathcal{W}$ in the following lemma.
\begin{lemma}\label{prvo:l_repVW}
Let $\x=\ra(P)$, $V=\mH^{1/2}P\mH^{'\dagger 1/2}$ and $W=\mH^{1/2}P_\perp\mH^{'\dagger 1/2}$
then
\begin{align}\label{prvo:e_piV}
V^*V &=P_{\ra(\mH' P)}\\
W^*W&= P_{\ra(\mH' P_\perp)} \label{prvo:e_piV3}\\
WV^*&=VW^*=0,\label{prvo:e_piV2}\\
\mathcal{W}&={\sf inv}(\mH^{1/2})\x,\label{prvo:e_inverseimage}
\end{align}
where $\mathcal{W}$ is from (\ref{drugo:eq:sefubacio}) and
$$
{\sf inv}(\mH^{1/2})\x=\{x : \mH^{1/2}x\in\x\}
$$
denotes the inverse image of the subspace $\x$ under the mapping $\mH^{1/2}$.\index{inverse image}
\end{lemma}
\begin{proof}
The relations (\ref{prvo:e_piV})--(\ref{prvo:e_piV2}) follow analogously as
in the proof of Theorem \ref{prvo1:l_lema13}. It only remains to prove
(\ref{prvo:e_inverseimage}).


We first show that ${\sf inv}(\mH^{1/2})\x\subset\mathcal{W}=\ra(W)^\perp$. Take any
$u\in{\sf inv}(\mH^{1/2})\x$, then
$$
\mH^{1/2}u=z\in\x.
$$
This implies
$$
0=(z,P_\perp\mH^{'\dagger1/2}v)=(u,\mH^{1/2}P_\perp\mH^{'\dagger1/2}v),\qquad v\in\H
$$
which proves $u\in\ra(W)^\perp=\mathcal{W}$.

The other inclusion follows in two steps.
Take $u\in\mathcal{W}$, then
$$
(u,\mH^{1/2}P_\perp\mH^{'\dagger1/2}v)=0,\qquad v\in\H.
$$
On the other hand, the subspace
$$\ra(P_\perp\mH^{'\dagger 1/2})^\perp=\ra(P_\perp P_{\ra(\mH')})^\perp\subset\d(\mH^{1/2})$$
is finite dimensional, so we conclude $u\in\d(\mH^{1/2})$.
Corollary \ref{prvo:corr:333} implies
$$
0=(\mH^{1/2}u, P_\perp P_{\ra(\mH')}v)=(\mH^{1/2}u, P_{\ra(\mH')} P_\perp v)=
(\mH^{1/2}u, P_\perp v),\qquad v\in\H,
$$
which proves $\mH^{1/2}u\in\x$. With
this conclusion we have established (\ref{prvo:e_inverseimage}).
\qquad\end{proof}

\bigskip
As a direct consequence of Corollary \ref{prvo:hoce_sef} and (\ref{prvo:e_inverseimage})
we obtain the following result.
\begin{corollary}\label{prvo:corollary:sinp}
Let $\x=\ra(P)$, $V=\mH^{1/2}P\mH^{'\dagger 1/2}$ and $W=\mH^{1/2}P_\perp\mH^{'\dagger 1/2}$
then $$\|P_\mathcal{V}P_{\mathcal{W}^\perp}\|\leq\|P_{\mathcal{V}^\perp}
P_{\mathcal{W}}\|,$$ so
\begin{equation}\label{prvo:e_dimenzijaVW}
\sin\Theta_p(\mH^{1/2}\x, {\sf inv}(\mH^{1/2})\x)=\|V^*W\|.
\end{equation}
\end{corollary}
It would be pleasing to use $\mH^{1/2\dagger}$
in the place of ${\sf inv} (\mH^{1/2})$. This
is only possible under additional restrictions on the subspace
$\ra(P)$. To get better feeling for the meaning of $\sinbf\Theta_p(\mH^{1/2}\x, {\sf inv}(\mH^{1/2})\x)$
consider the following example.
\begin{example}{\rm\label{example_1}
Take
$$
H=\left[\begin{matrix}1&1\\1&1\end{matrix}\right],\qquad\x=\left[\begin{array}{c}1\\0
\end{array}\right],
$$
then
$$
H'=\left[\begin{matrix}1&0\\0&1\end{matrix}\right]
$$
is, unlike $H$, a positive definite matrix.
Now,
$$
H^{1/2}=\left[\begin{matrix} \frac{1}{{\sqrt{2}}} & \frac{1}
   {{\sqrt{2}}} \cr \frac{1}{{\sqrt{2}}} & \frac{1}
   {{\sqrt{2}}}   \end{matrix}\right],\qquad
   H^{1/2\dagger}=\left[\begin{matrix} \frac{1}{2\,{\sqrt{2}}} & \frac{1}
   {2\,{\sqrt{2}}} \cr \frac{1}{2\,{\sqrt{2}}} &
    \frac{1}{2\,{\sqrt{2}}}  \end{matrix}\right],
$$
and we compute
$$
\ra(V)=\text{span}\{[\begin{matrix}1& 1\end{matrix}]^*\},\qquad
\ra(W)^\perp=\text{span}\{[\begin{matrix}-1& 1\end{matrix}]^*\}
$$
which proves that in this case $\sinbf\Theta_p(\ra(V), \ra(W)^\perp)=1$ and
$$
\ra(W)^\perp=\je(H)\not=\ra(\mH^{1/2\dagger}P).
$$
}\end{example}

Instead of advocating the use of the general formula (\ref{prvo:e_inverseimage})
 we will establish a ``compatibility condition'' under which we may use the generalized
inverse of $\mH^{1/2}$ to check the statement of the theorems.


The next result is a nonnegative analogue of Theorem \ref{prvo1:l_lema13}.
It will enable us to, in effect, ``deflate away" the kernel of the nonnegative form
$h$ and reduce the problem to the positive definite case.

\begin{theorem}\label{prvo1:t_dekompozicija}
Let the subspace $\x=\ra(P)\subset\q$ be given and let $h$ be a nonnegative form.
Assume $\sinbf\Theta_p(\mH^{1/2}\x,{\sf inv}(\mH^{1/2})\x)=\sinbf\Theta_p<1$ then
\begin{align}\label{prvo:e_prvizakjN}
(1-\sin\Theta_p)h'[u]
&\leq h[u]\leq(1+\sin\Theta_p)h'[u],\qquad u\in\q(h),\\
(1-\frac{\sin\Theta_p}{1-\sin\Theta_p})h[u]
&\leq h'[u]\leq(1+\frac{\sin\Theta_p}{1-\sin\Theta_p})
h[u],\qquad u\in\q(h).\label{prvo:e_drugizakjN}
\end{align}
\end{theorem}
\begin{proof}
The proof is similar to the proof of Theorem \ref{prvo1:l_lema13}. Let $h'$ and $\delta h$
be as in (\ref{prvo1:e_generickiH}). Set $\delta H_s$ to be the operator
defined by the form
$$
\delta h_s(x, y)=\delta h(\mH^{'\dagger 1/2}x, \mH^{'\dagger 1/2}y),\qquad x,y\in\H.
$$
%
The form $\delta h_s$ is closed and everywhere defined, so
$\delta H_s$ is a bounded operator.
We obviously have
$\je(\mH^{'\dagger 1/2})=\je(\mH')\subset\je(\delta H_s)$, so
$P_{\ra(\mH')}$ commutes with the operator $\delta H_s$.
With the use of Corollary \ref{prvo:corr:333}
one computes, analogously as in Theorem \ref{prvo1:l_lema13},
\begin{align*}
\delta h(\mH^{'\dagger 1/2}x, \mH^{'\dagger 1/2}y)&=h(P_\perp\mH^{'\dagger 1/2}x,
P\mH^{'\dagger 1/2}y)+ h(P\mH^{'\dagger 1/2}x, P_\perp\mH^{'\dagger 1/2}y)\\
&=(Wx, Vy)+(Vx, Wy),
\end{align*}
so
$$
\delta H_s=V^*W+W^*V.
$$
Since $\mH^{'1/2}\mH^{'\dagger 1/2}=P_{\ra(\mH')}$ we obtain
\begin{equation}\label{prvo:e_deltaJezgra}
\max_{u,v\in\ra(\mH')\cap\q}\frac{|\delta h(u, v)|}{\sqrt{h'[u]h'[v]}}=
\|\delta H_s\|=\|V^*W\|.
\end{equation}
Corollary \ref{prvo:corollary:sinp} implies that the assumption $\sinbf\Theta_p<1$,
in fact, reads
$$
\sin\Theta_p=\|V^*W\|<1.
$$
With this in hand, we have established
$$
(1-\sin\Theta_p)h'[u]
\leq h[u]\leq(1+\sin\Theta_p)h'[u],\qquad u\in\q(h),
$$
which implies $\je(\mH')=\je(\mH)$. The relation (\ref{prvo:e_drugizakjN})
follows by the same argument as the one used in Theorem \ref{prvo1:l_lema13}.
\qquad\end{proof}

The main insight into the structure of the operator $\mH'$,
gained from Theorem \ref{prvo1:t_dekompozicija}, is summed
up in the following corollary.

\begin{corollary}\label{c_redukcija}
Take a nonnegative form $h$ and a subspace $\x=\ra(P)\subset\q$.
If $\sinbf\Theta_p(\mH^{1/2}\x,{\sf inv}(\mH^{1/2})\x)<1$ then $\ra(\mH')=\ra(\mH)$.
\end{corollary}

Corollary \ref{c_redukcija} gives precise meaning to the statement
``deflate away''.
Set $\mathcal{R}=\ra(\mH)=\ra(\mH')$ and $\mathcal{N}=\je(\mH)=\je(\mH')$.
The projections $P_{\mathcal{N}}$ and $P$ commute,
so
$$
P_{\mathcal{N}\cap\ra(P)}=P_{\mathcal{N}}P,\qquad \widetilde{P}=P-P_{\mathcal{N}\cap\ra(P)}
$$
are both orthogonal projections. A direct calculation shows
$$
\widetilde{\x}:=\ra(\widetilde{P})=\ra(P)\ominus(\mathcal{N}\cap\ra(P))=\ra(\mH')\cap\ra(P)=\ra(\mH'P).
$$
The form
$$
\widetilde{h}(u,v)=h(P_{\mathcal{R}}u,P_{\mathcal{R}}v)
$$
is positive definite in $\mathcal{R}$ and $\ra(\widetilde{P})\subset\q(\widetilde{h})\cap\mathcal{R}$.
Now, apply the construction (\ref{prvo:e_def_dh})---(\ref{prvo1:e_generickiH}) to
the form $\widetilde{h}$ and the projection $\widetilde{P}$.
By $\widetilde{\mH}:\mathcal{R}\to\mathcal{R}$ denote
the operator defined by the form $\widetilde{h}$ in $\mathcal{R}$,
then $\ra(\widetilde{P})\subset\mathcal{R}$
and $$\widetilde{h}'(u,v)=h'(P_{\mathcal{R}}u,P_{\mathcal{R}}v).$$
We conclude that
$$
\sin\Theta(\widetilde{\mH}^{1/2}\widetilde{\x},\widetilde{\mH}^{-1/2}\widetilde{\x})
=\sin\Theta_p(\mH^{1/2}\x,{\sf inv}(\mH^{1/2})\x)<1
$$
and $\widetilde{h}$ and $\widetilde{P}$ satisfy the
assumptions of Theorem \ref{prvo1:l_lema13}. If we were to
``a priori'' assume $\ra(\mH')=\ra(\mH)$, then this
argument would give an alternative proof
of Theorem \ref{prvo1:t_dekompozicija}. ``Deflate away" means that we
assume we were given $\widetilde{h}$ and $\widetilde{P}$ as input.

\begin{remark}{\rm\label{rem:deflate_away}
Another consequence of Corollary \ref{c_redukcija} is that we can invoke Lemma \ref{prvo2:t_drmac}
to conclude that the constant $\frac{\sinbf\Theta_p}{1-\sinbf\Theta_p}$ (in (\ref{prvo:e_drugizakjN}))
cannot be sharpened.
Furthermore, Example \ref{example_1} shows that the assumption
$$\sin\Theta_p(\mH^{1/2}\x,{\sf inv}(\mH^{1/2})\x)<1$$ is a necessary requirement
to establish the inequalities (\ref{prvo:e_prvizakjN}) and (\ref{prvo:e_drugizakjN})
as well as to guarantee that $\ra(\mH)=\ra(\mH')$ (equivalently  $\je(\mH)=\je(\mH')$).
}\end{remark}

\subsubsection{Important special case}
The assumption that $P$ and $P_{\je(\mH)}$ commute and
Corollary \ref{prvo:corr:333} yield $\je(\mH)=\je(\mH')$ and $\ra(\mH)=\ra(\mH')$. This implies
\begin{equation}\label{prvo:eq:geninv}
{\sf inv}(\mH^{1/2})\x=\mH^{1/2\dagger}\x.
\end{equation}
The projections $P$ and $P_{\je(\mH)}$ certainly commute when $\je(\mH)\perp\ra(P)$
or when\footnote{The other situation when
$P$ and $P_{\je(\mH)}$ commute is when $\ra(P)\subset\je(\mH)$,
this situation is however trivial and we have tacitly left it out.}
$\je(\mH)\subset\ra(P)$. This discussion is summed up in the following corollary.
\begin{corollary}\label{prvo:cor:ginv}
Assume $P=XX^*$ and $P_{\je(\mH)}$ commute and let $$\sinbf\Theta_p(\mH^{1/2}X, \mH^{1/2\dagger}X)<1,$$
then
\begin{align}\label{prvo:e_prvizakjN:inv}
(1-\sin\Theta_p)h'[u]
&\leq h[u]\leq(1+\sin\Theta_p)h'[u],\qquad u\in\q(h)\\
(1-\frac{\sin\Theta_p}{1-\sin\Theta_p})h[u]
&\leq h'[u]\leq(1+\frac{\sin\Theta_p}{1-\sin\Theta_p})
h[u],\qquad u\in\q(h).\label{prvo:e_drugizakjN:inv}
\end{align}
\end{corollary}

\begin{remark}{\rm
To assess the restriction
that $P$ and $P_{\je(\mH)}$ should commute, consider the definition of the
relatively accurate approximation of the number $\lambda\in\R_{+}$. $\mu\in\R_{+}$
is relatively accurate approximation of $\lambda\in\R_{+}$,
if
\begin{enumerate}
  \item $\lambda=\mu$, when $\lambda=0$
  \item $\frac{|\lambda-\mu|}{\mu}<1$, when $\lambda\not= 0$.
\end{enumerate}
This implies that we can expect to compute ``relatively accurate''
Ritz value approximation of the spectrum of the nonnegative definite
operator $\mH$ only in the case when we have computed a basis for $\je(\mH)$,
cf. \cite{ArbenzDrmac}.}\end{remark}

Remark \ref{rem:deflate_away} implies that we may assume
that the condition of Corollary \ref{prvo:cor:ginv} were $\je(\mH)\perp\ra(P)$.
To compute the basis of the set ${\sf inv}(\mH^{1/2})\x$ we need to repeatedly solve the
equation
$$
\mH^{1/2}u=x_i,\qquad i=1,...,{\sf dim}(\x).
$$
The vectors $x_i$ are assumed to be a basis for $\x$. The
restriction that $\je(\mH)\perp\ra(P)$ amounts to nothing more
then to impose a compatibility condition on $x_i$ (e.g. think of
the Laplacian with Neumann boundary conditions).


\subsection{A first approximation estimate}
Theorem \ref{prvo:t_monotonost} and Lemma \ref{prvo1:l_reducira} yield
the first eigenvalue estimates.
The next theorem will give an eigenvalue estimate with the minimum of the
restrictions on the subspace $\ra(X)\subset\q$. Sharper bounds are possible
when we impose additional assumptions on $\ra(X)$. Even this (first order) estimate
will compare favorably with other higher order bounds that can be found in the literature,
cf. \cite{GruVes02}.

\begin{theorem}\label{prvo1:t_selekcija}
Let $0\leq h$ and let the $n$-dimensional subspace $\ra(P)\subset\q$, $P=XX^*$, be given.
Define
$$
\Xi=(\mH^{1/2}X)^*\mH^{1/2}X,\quad
~\Xi\in\C^{n\times n}
$$
and assume $\mu_n<\lambda_e(\mH)$.
Here, the Ritz values are numbered as in (\ref{prvo:e_uredjRitz}).
If $\ra(P)$ is such that $\sinbf\Theta_p<1$,
then there are $n$ eigenvalues of the operator $\mH$,
counting the eigenvalues according to their multiplicities,
such that
\begin{align}\label{prvo1:e_tvrdnjaza0}
\displaystyle {|\lambda_{i_j}-\mu_j|}%
&\leq{\mu_j}\sin\Theta_p,\qquad j=1,\ldots,n\\
\displaystyle {|\lambda_{i_j}-\mu_j|}%
&\leq{\lambda_{i_j}}\frac{\sin\Theta_p}{1-\sin\Theta_p},\qquad j=1,\ldots,n,
\label{prvo1:e_tvrdnjaza+eta}
\end{align}
where $i_{(\cdot)}:\N\to\N$ is a permutation.
\end{theorem}
\begin{proof}
Corollary \ref{c_redukcija} readily implies the conclusion (\ref{prvo1:e_tvrdnjaza0})
for the Ritz values $\mu_j=0$, $j=1,\ldots,{\sf dim} (\je(\Xi))$.
Therefore, we may safely assume that $h$ is a positive definite form.
Lemma \ref{prvo1:l_reducira} implies $\sigma_{ess}(\mH)=\sigma_{ess}(\mH')$,
so the assumption $\mu_n<\lambda_e(\mH)$ guarantees that $\mu_n$ is a discrete eigenvalue of $\mH'$.
Theorem \ref{prvo1:t_dekompozicija} established
\begin{align*}
(1-\sin\Theta_p)h'[u]
&\leq h[u]\leq(1+\sin\Theta_p)h'[u],\qquad u\in\q(h)\\
(1-\frac{\sin\Theta_p}{1-\sin\Theta_p})h[u]
&\leq h'[u]\leq(1+\frac{\sin\Theta_p}{1-\sin\Theta_p})
h[u],\qquad u\in\q(h).
\end{align*}
The conclusion follows directly from Theorem \ref{prvo:t_monotonost}.
\qquad\end{proof}

For the numerical evidence concerning the performance of the
estimate (\ref{prvo1:e_tvrdnjaza+eta}) see the numerical tests from \cite{GruVes02}.

\section{Localizing the approximated eigenvalues}\label{prvo:s:localization}

There is a multitude of ways to match the computed Ritz values to a part of the
spectrum of the operator $\mH$ of the same multiplicity. These approaches usually differ
with regard to the allowed amount of additional information about
the spectrum of the operator $\mH$. Here, we present
two possible answers to that problem.

Theorem \ref{prvo1:t_selekcija} can be interpreted as a first localization result.
It gives an estimate of the infimum
of  $$\max_{j=1,..., n}\frac{|\lambda_{i_j}-\mu_j|}{\mu_j}$$ over
all of the permutations $i_{(\cdot)}:\N\to\N$.
So, we would be correct in stating that the Ritz values are approximating
the eigenvalues of $\mH$ that are closest to $\sigma(\Xi)$.

Having only limited additional infirmation we got a limited answer. We know that
there is a collection of eigenvalues of operator $\mH$, having the joint multiplicity $n$,
that is being approximated by the Ritz values from the subspace $\ra(X)$. The
information we have on the location of those eigenvalues in the spectrum of $\mH$ is
only that they are the eigenvalues closest to computed Ritz values.

Only when we have additional information about the location of the part
of the spectrum \textit{we do not want to approximate},
we can guarantee that we are approximating the
part of the spectrum we are interested in. A best known example of
such estimates is Temple--Kato inequality.\index{Temple--Kato inequality}
Assume $\lambda_1<\lambda_2$
and let $u\in\d(\mH)$ be a unit vector such that $(u,\mH u)<\gamma\leq\lambda_2$ then
\begin{equation}\label{prvo:e_TempleKato}
(u,\mH u)\geq\lambda_1\geq(u,\mH u)-\frac{(\mH u, \mH u)-(u,\mH u)^2}{\gamma-(u,\mH u)}.
\end{equation}
For a proof see \cite{ReedSimonSvi}. The estimate (\ref{prvo:e_TempleKato})
is valid for a general selfadjoint operator $\mH$.
In the following we shall formulate another assumption with the same effect, namely
to separate the ``unwanted'' component of the spectrum from
the Ritz values. Our result, however, does not need the regularity
constraint $u\in\d(\mH)$. Moreover, we will
obtain sharp bounds for the matching cluster of eigenvalues.
In the last section of this chapter we will demonstrate that on some
examples our bound considerably outperforms the estimate (\ref{prvo:e_TempleKato}).

We now give a theorem that determines those eigenvalues of the operator $\mH$,
given by a symmetric form $h$, which are approximated by the Ritz values associated
with the test subspace $\ra(X)\subset\q$.
Before we proceed with the formulation of the theorem we state a well
known fact that given $0<\lambda,\mu$ and $\sin\Theta_p<1$ the relation
$$
\frac{|\lambda-\mu|}{\mu}\leq\sin\Theta_p<1
$$
implies the relation
\begin{equation}\label{prvo1:e_etadef}
\frac{|\lambda-\mu|}{\lambda}\leq\frac{\sin\Theta_p}{1-\sin\Theta_p}\leq 2\sin\Theta_p.
\end{equation}

\begin{theorem}\label{prvo1:t_mat_new}
Set $
\gamma_r=\min \big\{(\lambda_p-\mu_k){(\lambda_p+\mu_k)^{-1}}|k=1,...,n; \;p=n+1,...,\infty\big\}
$ and
$
\eta_{\Theta_p}=\sinbf\Theta_p(1-\sinbf\Theta_p)^{-1}.
$
Take a nonnegative form $h$ and the subspace $\ra(X)\subset\q$. Assume
$r={\rm dim}(\je(\mH))\leq n$, set $P=XX^*$ and
let $h'$ be as in (\ref{prvo:e_def_h'}). By
$
\mu_1\leq\cdots\leq\mu_n
$,
denote the eigenvalues of the matrix
$
\Xi=(\mH^{1/2}X)^*\mH^{1/2}X\in\C^{n\times n}.
$ If $\gamma_r\geq 0$ and
$
\eta_{\Theta_p}<\min\{\gamma_r,1\}
$
then
\begin{equation}\label{prvo1:e_ostatak}
|\lambda_i-\mu_i|\leq\mu_i\sin\Theta_p,\quad i=1,...,n.
\end{equation}
\end{theorem}

{\em Proof}.
The assumption $\eta_{\Theta_p}<\min\{\gamma_r,1\}$ and Theorem \ref{prvo1:t_dekompozicija}
imply $\je(\mH)\subset \ra(X)$. Also, by Theorem \ref{prvo1:t_dekompozicija}
we have $\je(\mH)=\je(\mH')$, so
we are allowed to ``deflate away" the kernel of $\mH$.
Therefore, set $P_1=P_{\ra(\mH'P)}$ and proceed as if $h$
were positive definite and $P=P_1$.

The rest of the proof is completely analogous to the proof of
\cite[Theorem 5.1]{GruVes02}. The only difference is that in the place of
$\eta=\sinbf\Theta_p/(1-\sinbf\Theta_p)$ from \cite[Theorem 5.1]{GruVes02}
one uses a sharper quantity $\sinbf\Theta_p$. \qquad\endproof

If we are provided with the information that
\begin{equation}\label{prvo1:e_ugnjezdeno}
\eta_{\Theta_p}=\frac{\sin\Theta_p}{1-\sin\Theta_p}<\gamma_c:=
\min\Big\{\min_{\substack{k=1,...,n\\p=1,...,q-1}}
\frac{\mu_k-\lambda_p}{\lambda_p+\mu_k}, \min_{\substack{k=1,...,n\\p=q+n,...,\infty}}
\frac{\lambda_p-\mu_k}{\lambda_p+\mu_k}, 1\Big\}
\end{equation}
then $\mu_1\leq\cdots\leq\mu_n$
approximate the ``inner" eigenvalues
$$
\lambda_{q}\leq \lambda_{q+2}\cdots \leq\lambda_{q+n-1}~.
$$
This statement is made precise in the following theorem.

\begin{theorem}\label{prvo1:t_uronjeni}
Take a nonnegative definite form $h$ and a subspace $\ra(X)\subset\q$.
By $\mu_1\leq\cdots\leq\mu_n$ denote the eigenvalues of the matrix
$
\Xi=(\mH^{1/2}X)^*\mH^{1/2}X\in\C^{n\times n}$.
If
$
\eta_{\Theta_p}<\gamma_c
$,
where $\gamma_c$ is as in (\ref{prvo1:e_ugnjezdeno}),
then $\ra(P)\subset\ra(\mH')$ and
$$
|\lambda_{i+q-1}-\mu_i|\leq{\mu_i}\sin\Theta_p,\quad i=1,...,n~.
$$
\end{theorem}
\begin{proof}
The assumption (\ref{prvo1:e_ugnjezdeno}) and Theorem \ref{prvo1:t_dekompozicija} and
Corollary \ref{c_redukcija} imply
$$\ra(\mH')=\ra(\mH)\qquad\text{and}\qquad\ra(P)\subset\ra(\mH).$$
The rest of the proof follows analogously as in the proof
of Theorem \ref{prvo1:t_mat_new}.
\qquad\end{proof}
\begin{remark}{\rm\label{prvo:remark:cardinality}
Theorems \ref{prvo1:t_mat_new} and \ref{prvo1:t_uronjeni} imply that the
spectrum of the operator $\mH$ can stably (sensibly) be divided in two
disjoint parts: the part that is being approximated by the
$\sigma(\Xi)$ and the rest of the spectrum. To understand
this statement assume that the conditions of
Theorem \ref{prvo1:t_mat_new} hold. In this case both of the ``block diagonal'' forms
\begin{align*}
h(u,v)&=h(E(\lambda_n)u, E(\lambda_n) v)+h(E(\lambda_n)_\perp u, E(\lambda_n)_\perp v)
\simeq\left[\begin{array}{cc}\Lambda&\\&\mathbf{\Lambda}_c\end{array}\right]\\
h'(u,v)&=h(P u,P v)+h(P_\perp u, P_\perp v)
\simeq\left[\begin{array}{cc}\Xi&\\&\mathbf{\Xi}_c\end{array}\right]
\end{align*}
have ``diagonal blocks'' with disjoint spectra. We have assumed
$\Lambda=\text{diag}(\lambda_1,\ldots,\lambda_n)$ and
$\Xi=\text{diag}(\mu_1,\ldots,\mu_n)$ and $\mathbf{\Xi}_c$ and
$\mathbf{\Lambda}_c$ were unbounded operators defined by the forms
$h'$ and $h$ in the spaces $\ra(P_\perp)$ and
$\ra(E(\lambda_n)_\perp)$. In fact, we will colloquially call $h'$
the \textit{block diagonal part of the operator $\mH$ with respect
to the subspace $\ra(P)$}\index{block diagonal part of an
operator}. We will use the notation $h_P$ to denote $h'$ in
situations when it is not clear with respect to which test space
$\ra(P)$ was this construction performed.
}\end{remark}
\section{Eigenvector approximation estimates}
For the computed Ritz values
$$
0,0,\ldots,0,\mu_{r+1},\mu_{r+2},\ldots,\mu_n
$$
Theorem \ref{prvo1:t_selekcija} guarantees the existence of the eigenvalues
$$
\lambda_{i_1}\leq\lambda_{i_2}\leq\cdots\leq\lambda_{i_n},
$$
that are being approximated by the Ritz values (provided $\sinbf \Theta_p<1$) in the sense of
\begin{align*}
\displaystyle {|\lambda_{i_j}-\mu_j|}&\leq{\mu_j}\sin\Theta_p,\qquad j=1,\ldots,n.
\end{align*}

Assume $v_1, \ldots,v_n$ are mutually orthogonal
eigenvectors that belong to
the eigenvalues $\lambda_{i_1}\leq\lambda_{i_2}\leq\cdots\leq\lambda_{i_n}$.
If the conditions of
Theorems \ref{prvo1:t_mat_new} and \ref{prvo1:t_uronjeni}
are satisfied Remark \ref{prvo:remark:cardinality} assures us that
$$
\text{span}\{v_1, ..., v_n\}=\ra(E(\{\lambda_{i_1}, \lambda_{i_2},
\cdots, \lambda_{i_n}\})).
$$
Here we have assumed that $\mH=\int\lambda~\text{d}E(\lambda)$. To ease the presentation we generically use
$$
\widehat{E}=E(\{\lambda_{i_1}, \lambda_{i_2}, \cdots, \lambda_{i_n}\})
$$
to denote the projection on the subspace that is selected by a result like Theorem \ref{prvo1:t_mat_new}.

The central role in the analysis of the eigenvector approximations will be played
by the following lemma.

\begin{lemma}\label{prvo1:l_mathves}
Let $h$ be a nonnegative form and let $\mH^\dagger$ be bounded. Take
 $\ra(P)\subset\q$ such that $\sinbf\Theta_p<1$ and define
$$
s(x,y)=\delta h(\mH^{\dagger 1/2}x,\mH^{'\dagger 1/2} y),\qquad x,y\in\H.
$$
The form $s$ defines a bounded operator $S$ and
\begin{align}\label{prvo1:e_prvaDH}
S&=\mH^{1/2}\mH^{'\dagger 1/2}-\overline{\mH^{\dagger 1/2}
\mH'^{1/2}}\\
|(x, Sy)|&=\abs{s(x,y)}\leq\frac{\sin\Theta_p}{\sqrt{1-\sin\Theta_p}}\norm{x}\norm{y},\quad x,y\in\H.
\label{prvo1:e_mathves}
\end{align}
\end{lemma}
\begin{proof}
The closed graph theorem implies that
the operator
$$
S=\mH^{1/2}\mH^{'\dagger 1/2}-\overline{\mH^{\dagger 1/2}\mH'^{1/2}}
$$
is bounded. Also, $\je(\mH)=\je(\mH')=\je(S)$ and $P_{\je(S)}$ commutes with $S$.
It is sufficient to prove the estimate for $x, y\in\ra(\mH)$.
The inequality (\ref{prvo:e_deltaJezgra})
gives
$$
|\delta h(\mH^{\dagger 1/2} x,\mH^{'\dagger 1/2} y)|\leq\sin\Theta_p\norm{ y}~
h'[\mH^{\dagger 1/2} x]^{1/2}.
$$
Analogously, (\ref{prvo:e_prvizakjN}) implies
\begin{equation}\label{prvo:e_ocjenaZaUpotr}
\|\mH'^{1/2}\mH^{\dagger 1/2}\|\leq\frac{1}{\sqrt{1-\sin\Theta_p}}~.
\end{equation}
Altogether, the estimate (\ref{prvo1:e_mathves}) follows.\qquad\end{proof}

The operator $S$ has the special structure. Assume $\mH' u=\mu u$ and $\mH v=\lambda v$,
then
\begin{align}
\nonumber (v, Su)&=\lambda^{1/2}(v, u)\mu^{1/2} -\lambda^{-1/2}(v, u)\mu^{1/2}\\
&=\frac{\lambda-\mu}{\sqrt{\lambda\mu}}(v,u)\label{prvo:eq_strukruraSa}~.
\end{align}
The equation (\ref{prvo:eq_strukruraSa}) introduces the distance function
$$
\frac{\lambda-\mu}{\sqrt{\lambda\mu}}
$$
that measures the distance between the Ritz values and the spectrum of the operator $\mH$.
This distance function will feature in the important role in the estimates that follow.
The next theorem extends the scope, as well as strengthens the
eigenvector estimate from \cite{GruVes02,MathVeselic98} and is even new in
the matrix case. It can be seen as the eigenvector companion result of
Theorem \ref{prvo1:t_selekcija}.

\begin{theorem}\label{prvo1:t_mathves}
Let $h$ be a nonnegative form, and  let $\ra(P)\subset\q$ be such that
it satisfies the assumptions of Theorem \ref{prvo1:t_selekcija}.
Let $u_1,\ldots,u_n$ the mutually orthogonal eigenvectors belonging
the eigenvalues $\mu_1,\ldots,\mu_n$ of $\mH'P$, then
there exist mutually orthogonal eigenvectors $v_1,\ldots,v_n$
of $\mH$, belonging to the eigenvalues
$\lambda_{i_1},\ldots\lambda_{i_n}$
and
\begin{equation}\label{prvo1:e_vektmathvesDH}
\norm{v_j-u_j}\leq\frac{\sqrt{2}~\sin\Theta_p}{\sqrt{1-\sin\Theta_p}}
\max_{k\ne j}\frac{\sqrt{\mu_j\lambda_{i_k}}}{|\lambda_{i_k}-\mu_j|}.
\end{equation}
The eigenvalues $\lambda_{i_j}$, $j=1,\ldots,r$ are numbered in the ascending order as
given by Theorem \ref{prvo1:t_selekcija}.
\end{theorem}
\begin{proof}
Assume $\mu_1=\cdots=\mu_r=0$. Corollary \ref{c_redukcija} implies that
$u_i\in\je(\mH)$ for $i=1,\ldots,r$ so we take
$$
v_i=u_i,\quad i=1,\cdots,r.
$$
For $v_j$, $j=r+1,\ldots n$ take any orthonormal set of eigenvectors belonging to
the eigenvalues $\lambda_{i_j}$, $j=r+1,\ldots,n$.
Since both $u_i$ and $v_i$, for $i=r+1,\ldots, n$ are perpendicular
to $\je(\mH)$ we may assume that $\mH$ is positive definite and
we are only given $u_i$, $i=r+1,\ldots n$ as test vectors.
Take $s$ from Lemma \ref{prvo1:l_mathves} and use (\ref{prvo1:e_prvaDH}) to compute
\begin{eqnarray*}
s(v_k,u_j)&=&\delta h(\mH^{-1/2}v_k,\mH^{'-1/2}u_j)=
\sk{v_k,\mH^{1/2}\mH^{'-1/2}u_j}-\sk{\mH^{'1/2}\mH^{-1/2}v_k,u_j}\\
&=&(\lambda_{i_k}^{1/2}\mu_j^{-1/2}-\lambda_{i_k}^{-1/2}\mu_j^{1/2})\sk{v_k,u_j}
\end{eqnarray*}
and
\begin{eqnarray*}
\sum_{k\ne j}|\sk{v_k,u_j}|^2&\leq&\max_{k\ne j}
\frac{\lambda_{i_k}\mu_j}{(\lambda_{i_k}-\mu_j)^2}\sum_{k \ne j}|s(v_k, u_j)|^2\leq\max_{k\ne j}
\frac{\lambda_{i_k}\mu_j}{(\lambda_{i_k}-\mu_j)^2}~\|S^* u_j\|^2\\
&\leq&\max_{k\ne j}
\frac{\lambda_{i_k}\mu_j}{(\lambda_{i_k}-\mu_j)^2}\frac{\sin^2\Theta_p}{1-\sin\Theta_p}.
\end{eqnarray*}
Scaling $v_j, u_j$ so that $\sk{v_j,u_j}\geq0$, we obtain
\begin{eqnarray*}
\norm{v_j-u_j}&=&\sqrt{2}\Big[1-\sk{v_j,u_j}\Big]^{1/2}=\sqrt{2}
\Big[1-\big[1-\sum_{k \ne j}|\sk{v_k,u_j}|^2\big]^{1/2}\Big]^{1/2}\\
&\leq&\sqrt{2}\Big[1-\big[1-\max_{k\ne j}
\frac{\lambda_{i_k}\mu_j}{(\lambda_{i_k}-\mu_j)^2}\frac{\sin^2\Theta_p}{1-\sin\Theta_p}\big]^{1/2}\Big]^{1/2}\\
&\leq&\frac{\sqrt{2}\;\sin\Theta_p}{\sqrt{1-\sin\Theta_p}}
\max_{k\ne j}\frac{\sqrt{\mu_j\lambda_{i_k}}}{|\lambda_{i_k}-\mu_j|}.
\end{eqnarray*}
This proves the lemma in the case in which $\sigma_e(\mH)=\emptyset$. In the general case we
use the formula
$$
\frac{\sqrt{\lambda_e\mu_j}}{\lambda_e-\mu_j}~\Big|\Big(E_{\mH^{1/2}}(\big[\sqrt{\lambda_e},\infty\big>)u_j,Su_j\Big)\Big|\geq
\Big|\Big(E_{\mH^{1/2}}(\big[\sqrt{\lambda_e},\infty\big>)u_j,u_j\Big)\Big|
$$
and analogous argument. \qquad\end{proof}

\section{A simple model problem}
\label{section_problem}

We will now present an application of our theory to the singularly perturbed Sturm-Liouville
eigenvalue problem (\ref{e_toplina_u}). Estimates (\ref{prvo:e_TempleKato}) and
\cite[Theorem 6.21(Kato--Temple eigenvector estimate)]{Chat83}
do not apply due
to overly stringent regularity assumptions on the test vector, cf. (\ref{e_toplina_u})--(\ref{drugo:e_testfunkcija_u}).

Consider the family of positive definite forms
\begin{equation}\label{e_toplina}
h_\eta(u,v)=h_b(u,v)+\eta^2 h_e(u,v)=\int_{0}^2 u'v'~dx +\eta^2\int_{1}^2 u'v'~dx,
\quad u,v\in H^1_0[0,2].
\end{equation}
By $\mH_\eta$ denote the positive definite operator which is defined by the form $h_\eta$ from  (\ref{e_toplina}).
We are interested in eigenvalues and eigenvectors of the operator $\mH_\eta$
for large $\eta$. Here, $H^1_0[0,2]$ denotes the first order Sobolev space with
zero trace on the boundary.

This is the eigenvalue
problem for the vibration of a highly inhomogeneous string. We are only considering an academic example where we
can efficiently compute all information we need. For more realistic
applications see \cite{GruPhd}.

\begin{figure}[t]
\begin{center}\includegraphics[width=6cm]{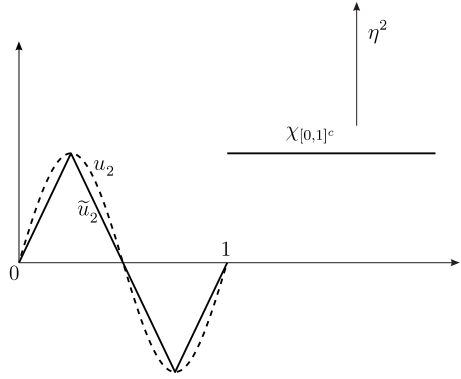}\end{center}
\caption[Square-well potential approximations]
{Various test functions for $\mH_\infty$ and $\mH_\eta$, $\eta$ large.}
\label{fig_boundeddomain}
\end{figure}

If we identify the functions from $H^1_0[0,\alpha]$, $\alpha>0$, with
their extension by zero to the whole of $[0,\beta]$ for $\beta\geq\alpha$, then we can
write
\begin{equation}\label{drugo:ulaganje}
H_0^1[0,\alpha]\subset H^1_0[0,\beta],\qquad 0<\alpha<\beta.
\end{equation}
Let $\chi_{[0,1]}$ be the characteristic function of the interval $[0,1]$ and
let $\chi_{[0,1]^c}=1-\chi_{[0,1]}$. Keeping (\ref{drugo:ulaganje}) in mind,
we conclude that
$$
\mH_\eta=-\partial_{x}(1+\eta^2\chi_{[0,1]^c})\partial_x,\qquad \d(\mH_\eta)=H^2[0,2]\cap H^1_0[0,2].
$$

It is known
that the forms $h_\eta$ converge to the form
$$
h_\infty(u,v)=\int^1_0 u'v'~dx,\quad u,v\in H^1_0\left[0,1\right]
$$
in the norm resolvent sense\footnote{More on the properties
of this convergence can be found in \cite{HempelPost,WeidmannScand84}.}.
Operators
$\mH_\eta$ and $\mH_\infty$ have discrete spectra and all the eigenvalues
are nondegenerate, cf. \cite{WeidmannSL}.
Since we will be considering the whole family of operators $\mH_\eta$,
additional notation will be introduced to ease the understanding. By
$$
\lambda_1^\eta<\cdots<\lambda_n^\eta<\cdots
$$
we denote the increasingly ordered eigenvalues of the operator $\mH_\eta$
and by
$$\lambda^\infty_1<\cdots<\lambda_n^\infty<\cdots
$$
the eigenvalues of the operator $\mH_\infty$.

The eigenpairs of
the operator $\mH_\infty$---which is defined in $L^2[0,1]$--- are \\$(n^2\pi^2, \sqrt{2}\sinbf(n \pi x))$, $n\in\N$.
The functions
\begin{equation}\label{drugo:e_testfunkcija}
u_n(x)=\begin{cases}\sqrt{2}\sin(n \pi x),&0\leq x\leq 1\\
0,& 1\leq x\end{cases},\qquad n\in\N
\end{equation}
are in $H^1_0[0,1]$ and also in $H^1_0[0,2]$.
 Therefore, they can be used as test functions
for an approximation of the eigenvalues of
$\mH_\eta$ ( for large $\eta$). Furthermore, according to
(\ref{intro:formula}) we obtain
$$
\sin^2\Theta_\eta(u_i):=\sin^2\Theta(\mH^{-1/2}_\eta u_i, \mH^{1/2}_\eta u_i)=
\frac{(u_i,\mH^{-1}_\eta u_i)-(u_i, \mH_\infty^{\dagger} u_i)}{(u_i,\mH^{-1}_\eta u_i)}.
$$

Let us now concentrate on the approximation of
the lowest eigenvalue. We compute
the Ritz value
$$
h_\eta(u_1,u_1)=\pi^2.
$$
When $\sinbf\Theta_\eta(u_1)<1$, Theorem \ref{prvo1:t_selekcija} guarantees the
existence of an eigenvalue
$\lambda^\eta_{i_1}$ such that
$$
\frac{|\lambda^\eta_{i_1}-\lambda^\infty_1|}{\lambda^\infty_1}\leq\sin\Theta_\eta(u_1)~.
$$
 A direct
computation shows that
\begin{align}
\nonumber (u_1,& \mH^{-1}_\eta u_1-\mH_\infty^{\dagger} u_1)\\
\nonumber &=\int^1_0\left[\int^x_0
2\left(\frac{y\,\left( 1 + (1+{\eta^2})\,\left( 1 - x \right)
\right) }{2 + {\eta^2}} - y\,\left( 1 - x \right)\right)
    \,\sin (\pi \,y)\sin (\pi \,x)~dy\right.\\
\nonumber &\left. +\int_x^1
2\,\left(\frac{x~\left( 1 + (1+{\eta^2})\,\left( 1 - y \right)
\right)}{2 + {\eta^2}} - x~\left( 1 - y \right)\right)
    \,\,\sin (\pi \,y)\sin (\pi \,x)~dy\right]~dx\\
    &=\frac{2}{(2+\eta^2)\pi^2}=O(\eta^{-2})\label{drugo:e_HighContrast}.
\end{align}

This establishes that $\sinbf\Theta_\eta(u_1)\to 0$, so Theorem \ref{prvo1:t_selekcija}
will be applicable for $\eta\geq 1$ such that
$$
\frac{(u_1,\mH^{-1}_\eta u_1)-(u_1, \mH_\infty^{\dagger} u_1)}{(u_1,\mH^{-1}_\eta u_1)}=
\frac{2}{4+\eta^2}<1.
$$
Furthermore, based on \cite{HempelPost} and \cite{WeidmannScand84} we conclude that the
assumptions of Theorem \ref{prvo1:t_mathves} must be satisfied for $\eta$ large.
We will now investigate this claim further.

The eigenvalues of the operator $\mH_\eta$ satisfy the equation
\begin{equation}\label{ritzII_eq}
\sqrt{1+\eta^2}\cot(\sqrt{\lambda^\eta})+\cot\big(\sqrt{\frac{\lambda^\eta}{1+\eta^2}}\big)=0.
\end{equation}
and the nonnormalized eigenvectors are
$$
\widehat{v}^\eta_i(x)=\begin{cases}
\sin(\sqrt{\lambda^\eta_i}x),&0\leq x\leq 1\\
\displaystyle\frac{\sin(\sqrt{\lambda^\eta_i})}{\sin\big(\sqrt{\frac{\lambda^\eta_i}{1+\eta^2}}\big)}
~\sin\big(\sqrt{\frac{\lambda^\eta_i}{1+\eta^2}}x\big),&1\leq x\end{cases}
$$
Set $v^\eta_i=\|\widehat{v}^\eta_i\|^{-1}~\widehat{v}^\eta_i$ then Theorems \ref{prvo1:t_mat_new} and \ref{prvo1:t_mathves}
imply
$$
\sin\angle(v^\eta_1,u_1)=\|v^\eta_1-u_1\|\leq\frac{\pi\sqrt{\lambda^\eta_2}}{\lambda^\eta_2-\pi^2}~\frac{2}{\sqrt{4+\eta^2-\sqrt{8+2\eta^2}}}.
$$
From (\ref{drugo:e_HighContrast}) we establish the uniform estimate
$$
\|v^\eta_1-u_1\|\leq
~\frac{1.333334}{\sqrt{4+\eta^2-\sqrt{8+2\eta^2}}},\qquad \eta\geq 2.
$$
This illustrate a way to obtain rigorous eigenvector estimates. First, we have localized the approximated
eigenvalue by an application of Theorem \ref{prvo1:t_mat_new}. This has selected the approximated eigenvector.
Theorem \ref{prvo1:t_mathves} then yields an accuracy of that approximation.

Let us note that
\begin{align*}
h_\infty(u_n,u_n)&=h_\eta(u_n,u_n)=n^2\pi^2\\
\sin\Theta_\eta(u_n)&=\frac{(u_n,\mH^{-1}_\eta u_n)-(u_n, \mH_\infty^{\dagger} u_n)}{(u_n,\mH^{-1}_\eta u_n)}=
\frac{2}{4+\eta^2}.
\end{align*}
This implies that we can get estimates for all
$\lambda^\eta_i$ and $v^\eta_i$ by an analogous procedure. In establishing
the convergence results for higher eigenvalues and eigenvectors it was important that we a priory
new that all $\lambda^\eta$ were nondegenerate. Our theory
has successfully been applied to similar singularly perturbed operators
which were defined in $L^2(\Omega)$, $\Omega\subset\R^n$, see \cite{GruPhd}.
For those operators such a claim does not hold. There it is important to generalize the
subspace results from \cite{GruVes02}
as well as to obtain higher order estimates (in $\sinbf\Theta_\eta$)
for eigenvalues. These results were obtained in the PhD. Thesis \cite{GruPhd}
and will be reported elsewhere.
\section{Conclusion}
A method to compute an estimate of the accuracy of the subspace
approximation method is presented. It can also be used to obtain
accurate lower estimates of a desired group of eigenvalues. The
bounds have to be viewed as a combination of the Ritz value bound,
which gives an existence of the matching of the Ritz values and
eigenvalues, and the subspace bound, which describes the nature
of that matching.

The case study that was just
performed can be described as leading to a ``pseudo spectral'' method.
We have used the completely solvable (``well behaved'') operator
$$
(\mH_\infty^{1/2}u, \mH^{1/2}_\infty v)
=h_\infty(u,v)=\int^1_0 u'v'~dx,\quad u,v\in H^1_0\left[0,1\right],
$$
to analyze the singularly perturbed operator $\mH_\eta$.
Since the eigenvalue problem for the operator $\mH_\infty$ was completely solvable,
we have used the eigenfunctions of the operator $\mH_\infty$ to define
a test space for the operator $\mH_\eta$. Analogously, we could have used
other test functions from $H^1_0[0,1]$ to analyze the operator $\mH_\eta$.
For instance, assume we have used the linear finite elements
to compute an approximation $\widetilde{u}_i$ of the
function $u_i$, see Figure \ref{fig_boundeddomain}.
Theorem \ref{prvo1:t_selekcija} can be invoked  if we find
a way to estimate
$\sinbf\Theta(\mH^{-1/2}_\eta \widetilde{u}_i,\mH^{1/2}_\eta \widetilde{u}_i)$.
The study of singularly perturbed eigenvalue problems and finite
element spectral approximations has been performed in \cite{GruPhd}.
The results will be presented in subsequent reports.
\section*{Acknowledgement}The author would like to
thank Prof. Dr. Kre\v{s}imir Veseli\'{c}, Hagen for helpful discussions and
support during the research and the preparation of this manuscript. The author
also thanks Dr. Josip Tamba\v{c}a, Zagreb for stimulating discussions.
Partial support of the grant Nr. 0037122 of Ministry of
Science, Education and Sport, Croatia is also acknowledged.

\bibliographystyle{abbrv}
\def\cprime{$'$} \def\cprime{$'$} \def\cprime{$'$}

\end{document}